\documentclass[draft,leqno,12pt]{article}
\usepackage{amsthm, amsmath,amssymb}
\usepackage{indentfirst}

\theoremstyle{plain}
\newtheorem{thm}{Theorem}[section]
\newtheorem{cor}{Corollary}[section]

\theoremstyle{definition}
\newtheorem*{rem}{Remark}

\theoremstyle{remark}
\newtheorem{case}{Case}

\def\Cset{\mathbb{C}}

\title{SEMILOCAL CONVERGENCE OF TWO ITERATIVE METHODS FOR SIMULTANEOUS COMPUTATION
OF POLYNOMIAL ZEROS\thanks{This paper is published in: C. R. Acad. Bulg. Sci 59 (2006), No 7, 705--712.}}
\author{Petko~D.~Proinov}
\date{}

\begin{document}
\maketitle

\begin{abstract}
In this paper we study some iterative methods for simultaneous approximation of polynomial zeros. We give new semilocal convergence theorems with error bounds for Ehrlich's and Nourein's iterations.
Our theorems generalize and improve recent results of
\textsc{Zheng} and \textsc{Huang} [J. Comput. Math. 18 (2000), 113--122],
\textsc{Petkovi\'c} and \textsc{Herceg} [J. Comput. Appl. Math. 136 (2001), 283--307] and
\textsc{Nedi\'c} [Novi Sad J. Math. 31 (2001), 103--111].
We also present a new sufficient condition for simple zeros of a polynomial.

\textbf{Key words:} polynomial zeros, simultaneous methods, semilocal convergence,
Ehrlich method, B\"orsch-Supan method, Nourein method

\textbf{2000 Mathematics Subject Classification:} 65H05, 12Y05, 26C10, 30C15
\end{abstract}

\section{Introduction and notations}

Let $f$  be a monic complex polynomial of degree $n \ge 2$. A point $\xi$ in
${\Cset}^n$ is said to be a root-vector of $f$ if its components are exactly the
zeros of $f$ with their multiplicities. There are a lot of iterations for
simultaneous computation of all zeros of $f$ (see the monographs of \textsc{Sendov,
Andreev, Kjurkchiev} \cite{SAK94}, \textsc{Petkovi\'c, Herceg, Ilic} \cite{PHI97} and
\textsc{Kyurkchiev} \cite{Kyu98}). The famous one is \textsc{Weierstrass'} iteration
\cite{Wei03}
\begin{equation}  \label{Wei}
z^{k + 1}  = z^k - W(z^k ), \qquad k = 0,1,2, \ldots,
\end{equation}
where the operator $W$ in ${\Cset}^n$ is defined by
$W(z) = (W_1(z), \ldots, W_n(z))$ with
\[
W_i(z) = \frac{f(z_i)}{\prod \limits_{j \ne i}{(z_i - z_j)}}
\qquad (i = 1, 2, \cdots, n).
\]
Another iteration for simultaneous finding all zeros of $f$ is
\textsc{Ehrlich'}s iteration \cite{Ehr67}
\begin{equation}  \label{Mae}
z^{k + 1}  = F(z^k ), \qquad k = 0,1,2, \ldots,
\end{equation}
where the operator $F$ in ${\Cset}^n$ is defined by
$F(z) = (F_1(z), \ldots, F_n(z))$ with
\begin{equation} \label{MIF}
F_i(z) = z_i -
\frac{f(z_i)}{f'(z_i) - f(z_i) \sum\limits_{j \ne i}{1/(z_i  - z_j)}} \, .
\end{equation}
\textsc{Werner} \cite{Wer82} has proved that the iteration function $F$ can also be
written in the form
\begin{equation} \label{BSIF}
F_i(z) = z_i  - \frac{W_i(z)}
{1 + \sum\limits_{j \ne i} {W_j(z) / (z_i - z_j )}}.
\end{equation}
Ehrlich's method \eqref{Mae} with the iteration function $F$ defined by \eqref{BSIF}
instead of \eqref{MIF} is  known as B\"orsch-Supan's method since in such form it was
proposed for the first time by \textsc{B\"orsch-Supan} \cite{BS70}. The following
iteration is due to \textsc{Nourein} \cite{Nou77}
\begin{equation}  \label{Nou}
z^{k + 1}  = G(z^k ), \qquad k = 0,1,2, \ldots,
\end{equation}
where the operator $G$ in ${\Cset}^n$ is defined by
$G(z) = (G_1(z), \ldots, G_n(z))$ with
\begin{equation} \label{NIF}
G_i(z) = z_i - \frac{W_i(z)}
{1 + \sum\limits_{j \ne i} {W_j(z) / (z_i - z_j - W_i(z))}}.
\end{equation}
Nourein's method is also known as B\"orsch-Supan method with Weierstrass' correction.

Since 1996, some authors \cite{Pet96,PH96,PI97,PH97,PHI98,ZH00,PH01,Ned01}
have obtained semilocal convergence theorems for Ehrlich's and Nourein's methods from data at one point. The best results on Ehrlich' method are due to \textsc{Zheng} and \textsc{Huang} \cite{ZH00} and \textsc{Petkovi\'c} and
\textsc{Herceg} \cite{PH01}. The best results on Nourein's method are due to \textsc{Zheng} and \textsc{Huang} \cite{ZH00} and
\textsc{Nedi\'c} \cite{Ned01}.

In this paper, we present new semilocal convergence theorems for Ehrlich's and Nourein's iterations which generalize and improve all previous results in this area. We also present a new sufficient condition for simple zeros of a polynomial. The main results of the paper (Theorems \ref{thm:loc}, \ref{thm:Ehr} and \ref{thm:Nou}) will be proved elsewhere.

Throughout the paper the norm $\|.\|_p \, \, (1 \le p \le \infty)$
in ${\Cset}^n$
is defined as usual, i.e.
$\|z\|_p = \left( \sum\nolimits_{i = 1}^n {|z_i|^p} \right)^{1/p}$.
For a given point $z$ in ${\Cset}^n$ we define
\begin{equation}
d(z) = (d_1(z), \ldots, d_n(z)) \quad\text{and}\quad
\delta(z) = \min \{d_1(z), \ldots, d_n(z) \} ,
\end{equation}
where
$d_i(z) = \min \{ |z_i  - z_j | : j \ne i \}$
for $i = 1,2, \ldots, n$.
For a given point $z$ in ${\Cset}^n$ with distinct components we use the notations
\begin{equation}
\frac{W(z)}{d(z)} =
\left( {\frac{W_1(z)}{d_1(z)}, \ldots, \frac{W_n(z)}{d_n(z)}} \right)
\quad\text{and}\quad E(z) = \left\|  \frac{W(z)}{d(z)} \right\|_p .
\end{equation}

\section{Localization of simple polynomial zeros}

The following theorem is an improvement of a result of \textsc{Zheng} \cite{Zhe87}.

\begin{thm} \label{thm:loc} 
Let $f$ be a monic polynomial of degree $n \ge 2$.
Suppose there exist $1 \le p \le \infty$ and a point $z^0$ in ${\Cset}^n$ with distinct components such that
\begin{equation} \label{simple}
E(z^0) < 1/2^{1/q} \quad\text{and}\quad \phi(E(z^0)) < 1 .
\end{equation}
where $q$ is defined by $1/p + 1/q = 1$ and $\phi$ is
a real function defined on $[0,1/2^{1/q})$ by
\begin{equation} \label{eq:phiW}
\phi(x) = \frac{(n - 1)^{1/q} x}{(1 - x)(1 - 2^{1/q} x)}
\left (1 + \frac{x}{(n - 1)^{1/p} (1 - 2^{1/q} x)} \right )^{n - 1} .
\end{equation}
Then the following statements hold true.
\begin{enumerate}
\item The polynomial $f$ has only simple zeros.
\item The closed disks
\[
D_i = \{z \in \Cset : |z-(z_i^0 - W_i(z^0))| \le C \, |W_i(z^0)| \} ,
\quad i = 1,2, \ldots, n ,
\]
where $C= \theta \lambda / (1 - \theta \lambda^2)$,
$\lambda = \phi(E(z^0))$, $\theta  = 1-2^{1/q} E(z^0)$,
are mutually disjoint and each of them contains exactly one zero of $f$.
\end{enumerate}
\end{thm}

Note that under the conditions \eqref{simple} Weierstrass' method \eqref{Wei} is convergent with the second order of convergence
(see \textsc{Proinov} \cite{Pro06}).

\section{Semilocal convergence of Ehrlich's method}

\begin{thm} \label{thm:Ehr}  
Let $f$ be a monic polynomial of degree $n \ge 2$,
$1 \le p \le \infty$ and
$1/p + 1/q = 1$.
Define the real function
\begin{equation} \label{phiM}
\phi(x) = \frac{a x^2 }
{(1 - (a + 1) x) (1 - (a + b) x)}
\left (1 + \frac{1}{n-1} .
\frac{ax}{1 - (a + b) x}
\right )^{n - 1} ,
\end{equation}
where $a = (n - 1)^{1/q}$ and $b = 2^{1/q}$.
Suppose that $z^0$ is an initial point in ${\Cset}^n$ with distinct components  satisfying
\[
E(z^0) < 1 / (a + b) \quad\text{and}\quad \phi(E(z^0)) \le 1 .
\]
Then the following statements hold true.
\begin{enumerate}
\item Ehrlich's iterative sequence \eqref{Mae} is well-defined and convergent to a root-vector $\xi $ of $f$. Moreover, if $\phi (E(z^0 )) < 1$, then the order of convergence is three.
\item For each $k \ge 1$ we have the following a priori error estimate
\begin{equation}  \label{priM}
\left \| z^k - \xi \right \|_p \le A_k \,
\frac{{\theta} ^k {\lambda }^{(3^k - 1) / 2}}
{1 - \theta {\lambda} ^{3^k } } \left \| W(z^0) \right \|_p ,
\end{equation}
where $\lambda  = \phi(E(z^0))$, $\theta  = \psi(E(z^0))$,
$A_k = \mu(E(z^0) \lambda^{(3^k - 1) / 2})$
and the real  functions $\psi$ and $\mu$ are defined by
\[
\psi(x) = \frac{1 - (a + b) x}{1 - a x}
\quad\text{and}\quad \mu(x) = \frac{1}{1 - a x}.
\]
\item For all $k \ge 0$ we have the following a posteriori error estimate
\begin{equation}  \label{posM}
\left \| z^k - \xi \right \|_p \le \frac{{\mu}_k}
{1 - {\theta}_k {\lambda}_k} \left \| W(z^k) \right \|_p ,
\end{equation}
where $\lambda_k  = \phi(E(z^k))$, $\theta_k  = \psi(E(z^k))$ and
$\mu_k = \mu(E(z^k)$.
\end{enumerate}
\end{thm}

Setting $p = \infty$ in Theorem~\ref{thm:Ehr} we obtain the following corollary.

\begin{cor} \label{cor:EhrInf1}  
Let $f$ be a monic polynomial of degree $n \ge 2$ and let $C$ be a real number
satisfying
\begin{equation} \label{MaeC}
0 \le C < \frac{1}{n+1}
\quad\text{and}\quad
\frac{(n - 1)C^2}{(1 - n C)(1 - (n+1) C)}
\left(\frac{1 - n C}{1 - (n +1)C} \right)^{n - 1} < 1.
\end{equation}
Suppose that $z^0$ is an initial point in $\Cset^n $ satisfying
\begin{equation} \label{iniinf}
\left\| \frac{W(z^0)}{d(z^0)} \right\|_{\infty} \le C.
\end{equation}
Then Ehrlich's method \eqref{Mae} is convergent to a root-vector
$\xi $ of $f$ with the third order of convergence. Moreover, the error bounds \eqref{priM} and \eqref{posM} hold with $p=\infty$.
\end{cor}

Corollary \ref{cor:EhrInf1} improves Theorem 4.1 of \textsc{Petkovi\'c} and \textsc{Herceg} \cite{PH01}. Note that Petkovic and Herceg have proved only linear convergence of Ehrlich's method under the stronger condition
\(
\left\| W(z^0) \right\|_{\infty} \le C \; \delta(z^0)
\)
with $C$ satisfying \eqref{MaeC} as well as
\begin{equation} \label{MaePIC1}
\beta : = \frac{(n - 1)C^2 (1 + (n-1) C) }{(1 - n C)(1 - (n-1) C)}
\left(\frac{1 - n C}{1 - (n +1)C} \right)^{n - 1} < 1
\quad\text{and}
\end{equation}
\begin{equation}  \label{MaePIC2}
g(\beta) < \frac{1 - (n-1) C}{2C}
\quad\text{where}\quad
g(x) = \left \{ \begin{array}{ll}
{1 + 2 x} & \mbox{ for $0 < x \le 1/2$} \, ;\\
{1/(1-x)} & \mbox{ for $1/2 < x < 1$}.
\end{array} \right.
\end{equation}
Corollary \ref{cor:EhrInf1} shows that Petkovic-Herceg's assumptions \eqref{MaePIC1} and \eqref{MaePIC2} can be omitted.

The following result improves and generalizes Theorem~1 of
\textsc{Zheng} and \textsc{Huang} \cite{ZH00} as well as previous results
\cite{Pet96,PH96,PI97,PH97}.

\begin{cor} \label{cor:EhrInfLp}  
Let $f$ be a monic polynomial of degree $n \ge 3$,
$1 \le p \le \infty$ and $1/p + 1/q = 1$.
Suppose $z^0$ is an initial point in
$\Cset^n $ satisfying
\begin{equation}  \label{iniMNLp}
\left \| \frac{W(z^0)}{d(z^0)} \right \|_p \le
\frac{1}{2(n - 1)^{1/q} + 2} \, .
\end{equation}
Then $f$ has only simple zeros and Ehrlich' method \eqref{Mae} converges to
a root-vector $\xi$ of $f$ with the third order of convergence.
Moreover, the error estimates \eqref{priM} and \eqref{posM} hold.
\end{cor}

\begin{proof}
From Theorem~\ref{thm:loc} and the proof of Corollary 2 of \cite{Pro06} we conclude
that $f$ has only simple zeros. It is easy to compute that
\[
\phi(R) = \frac{a}{(a + 1)(a + 2 - b)}
\left (1+ \frac{1}{n-1} . \frac{a}{a + 2 - b} \right)^{n-1}.
\]
for $R=1/(2a+2)$, where $a$, $b$ and $\phi$ are defined as in
Theorem~\ref{thm:Ehr}.
By Theorem~\ref{thm:Ehr} it suffices to prove that $\phi(R) \le 1$.
If $a > e-1$, then
$\phi(R) < e / (a+1) < 1$.
If
$a \le e-1$, then
\(
\phi(R) < (e-1)/(a-b+2) \le (e-1)/2 < 1
\)
which completes the proof.
\end{proof}

\begin{rem}
In Corollary~\ref{cor:EhrInfLp} we consider an initial condition of the type
\[
E(z^0) \le \frac{1}{A(n-1)^{1/q}+B}.
\]
Note that such initial conditions can be obtained for every $A>1$. For example,
if $A>1$, then one can take
\[
B=2^{1/q}+\frac{1}{4(A-1)} \exp{\frac{1}{A-1}} \, .
\]
\end{rem}

\begin{cor} \label{cor:EhrInf2}  
Let $f$ be a monic polynomial of degree $n \ge 2$.
Suppose that $z^0$ is an initial point in $\Cset^n $ satisfying
\[
\left \| \frac{W(z^0)}{d(z^0)} \right \|_{\infty} \le \frac{1}{1.5 n + 1.8} \, .
\]
Then Ehrlich's method \eqref{Mae} is convergent to a root-vector $\xi$ of $f$ with the
third order of convergence. Moreover, the error estimates \eqref{priM} and \eqref{posM}
hold with $p=\infty$.
\end{cor}

\begin{proof}
Define $\phi_n = \phi(1/(1.5 n + 1.8))$, where $\phi$ is defined by \eqref{phiM} with
$p=\infty$. The sequence $(\phi_n )$ is increasing for $2 \le n \le 10$ and
decreasing for $n \ge 10$. Hence $\phi_n  \le \phi_{10} < 1$ for $n \ge 2$. Now the
conclusion follows from Theorem~\ref{thm:Ehr}.
\end{proof}

Corollary~\ref{cor:EhrInf2} improves Theorem~4.2 of \textsc{Petkovi\'c} and
\textsc{Herceg} \cite{PH01}. Note that the authors of this work have proved that
Ehrlich's method is convergent under the condition
\(
\| W(z^0) \|_{\infty} \le C(n) \, \delta(z^0),
\)
where
\[
C(n) = \left \{ \begin{array}{ll}
{{1 / (n + 4.5)}} & \text{for } \, n = 3, 4 ;\\
{{1 / (1.545 \, n + 5)}} & \text{for } \, n \ge 5 .
\end{array} \right.
\]

\begin{cor} \label{cor:EhrL1}  
Let $f$ be a monic polynomial of degree $n \ge 2$.
Suppose that $z^0$ is an initial point in $\Cset^n $ satisfying
\[
\left \| \frac{W(z^0)}{d(z^0)} \right \|_{1} \le R = 0.2922 \ldots .
\]
where $R$  is the unique positive solution of the equation
\begin{equation} \label{eqMRn1}
{\left( \frac{x}{1-2x} \right)}^2 \exp \frac{x}{1-2x} = 1.
\end{equation}
Then $f$ has only simple zeros and Ehrlich' method \eqref{Mae} converges to a root-vector
$\xi$ of $f$ with the third order of convergence and the estimates \eqref{priM} and
\eqref{posM} hold with $p=1$.
\end{cor}

\begin{proof}
It follows from Theorem~\ref{thm:loc} and the proof of Corollary 3 of \cite{Pro06} that
$f$ has only simple zeros. It is easy to show that
$\phi(x) < g(x)$ for $0 < x < 1/2$, where $\phi(x)$ is defined by \eqref{phiM} with $p=1$
and $g(x)$ denotes the left-hand side of the equation \eqref{eqMRn1}. Therefore, $\phi(R)< 1$
which according to Theorem~3.1 completes the proof .
\end{proof}

\section{Semilocal convergence of Nourein's method}

\begin{thm} \label{thm:Nou}  
Let $f$ be a monic polynomial of degree $n \ge 2$,
$1 \le p \le \infty$ and $1/p + 1/q = 1$. Define the real function
\begin{equation}  \label{phiN}
\phi(x) = \frac{a^2 x^3}{(1 - m x + b x^2)
(1 - (a + 2) x + x^2)}
\left (1 + \frac{1}{n-1} .
\frac{a(x - x^2)}{1 - m x + b x^2} \right )^{n - 1} ,
\end{equation}
where $a = (n - 1)^{1/q}$, $b = 2^{1/q}$ and $m=a+b+1$.
Suppose that $z^0$ is an initial point in ${\Cset}^n$
with distinct components satisfying
\[
E(z^0) < 2 / (m + \sqrt{m^2 - 4 b} \, )
\quad\text{and}\quad \phi(E(z^0)) \le 1 .
\]
Then the following statements hold true.
\begin{enumerate}
\item Nourein's iterative sequence \eqref{Nou} is well-defined and convergent to a
root-vector $\xi$  of $f$. Moreover, if $\phi (E(z^0 )) < 1$, then the order of
convergence is four.
\item For each $k \ge 1$ we have the following a priori error estimate
\begin{equation}  \label{priN}
\left \| z^k  - \xi \right \|_p  \le
\mu_k \, \frac{{\theta} ^k {\lambda }^{(4^k - 1) / 3}}
{1 - \theta {\lambda} ^{4^k } } \left \| W(z^0) \right \|_p ,
\end{equation}
where $\lambda  = \phi(E(z^0))$, $\theta  = \psi(E(z^0))$,
$\mu_k = \mu(E(z^0) \lambda^{(4^k - 1) / 3})$ and the real  functions $\psi$ and $\mu$ are
defined by
\[
\psi(x) = \frac{1 - m x + b x^2}{1 - (a +1) x}
\quad\text{and}\quad
\mu(x) = \frac{1 - x}{1 - (a +1) x}.
\]
\item For all $k \ge 0$ we have the following a posteriori error estimate
\begin{equation}  \label{posN}
\left \| z^k - \xi \right \|_p \le \frac{{\mu}_k}
{1 - {\theta}_k {\lambda}_k} \left \| W(z^k) \right \|_p ,
\end{equation}
where $\lambda_k  = \phi(E(z^k))$, $\theta_k  = \psi(E(z^k))$ and
$\mu_k = \mu(E(z^k)$.
\end{enumerate}
\end{thm}

Setting $p = \infty$ in Theorem~\ref{thm:Nou} we obtain the following corollary.

\begin{cor} \label{cor:NouInf1}  
Let $f$ be a monic polynomial of degree $n \ge 2$ and let $C$ be a real number
such that $0 \le C < 2 / (n + 2 + \sqrt{n^2 + 4 n-8} \, )$ and
\begin{equation} \label{NouC}
\frac{(n - 1)^2 C^3}{(1 - (n+1) C+ C^2)(1 - (n +2) C+ 2 C^2)}
\left(\frac{1 - (n+1) C+ C^2}{1 - (n +2) C+ 2 C^2} \right)^{n - 1} < 1.
\end{equation}
Suppose that $z^0$ is an initial point in $\Cset^n $ satisfying \eqref{iniinf}.
Then Nourein's method \eqref{Nou} is convergent to a root-vector
$\xi $ of $f$ with the fourth order of convergence. Moreover, the error bounds \eqref{priN}
and \eqref{posN} hold with $p=\infty$.
\end{cor}

Corollary~\ref{cor:NouInf1} improves Theorem 2 of \textsc{Nedi\'c} \cite{Ned01}.
Nedi\'c has proved only linear convergence of Nourein's method under the condition
$\| W(z^0) \|_{\infty} \le C \; \delta(z^0)$ with $C$ satisfying
\[
0 < C < \frac{2}{n + 4 + \sqrt{n^2 + 8 n}} \quad\text{and}\quad \beta < \frac{1 - n
C}{1 + (n-2) C} \, ,
\]
where $\beta$ denotes the left-hand side of \eqref{NouC}.

The following corollary improves and generalizes Theorem~2 of
\textsc{Zheng} and \textsc{Huang} \cite{ZH00} as well as previous results
\cite{PH96,PH97,PHI98}.

\begin{cor} \label{cor:NouLp}  
Under the assumptions of Corollary~\ref{cor:EhrInfLp} $f$ has only simple zeros and
Nourein's method \eqref{Nou} converges to a root-vector $\xi$ of $f$ with the fourth order
of convergence.
Moreover, we have the error estimates \eqref{priN} and \eqref{posN}.
\end{cor}

\begin{proof}
Corollary~\ref{cor:EhrInfLp} implies that $f$ has only simple zeros. It it is easy to
compute that
\[
\phi(R) = \frac{2a^2(a+1)}{(2(a+1)^2 -(2a+1)b)(2a^2+2a+1))}
\left (1+ \frac{1}{n-1} . \frac{a(2a+1)}{2(a+1)^2 -(2a+1)b} \right)^{n-1}
\]
for $R=1/(2a+2)$, where $a$, $b$ and $\phi$ are defined as in
Theorem~\ref{thm:Nou}.
Taking into account that $b \le 2$ we get
\begin{equation} \label{b1}
\phi(R) \le \frac{2a^2(a+1)}{(2(a+1)^2 -(2a+1)b)(2a^2+2a+1))}
\left (1+ \frac{1}{n-1} \frac{2a+1}{2a} \right)^{n-1}
\end{equation}
and
\begin{equation} \label{b11}
\phi(R) \le \frac{a+1}{2a^2+2a+1}
\left (1+ \frac{1}{n-1} . \frac{2a+1}{2a} \right)^{n-1} .
\end{equation}
By Theorem~\ref{thm:Nou} it suffices to prove that $\phi(R) \le 1$. We shall consider
two cases.

\begin{case}
Suppose $a \ge 1.8$. It follows from \eqref{b11} that
\begin{equation}
\phi(R) < g(a) = \frac{a+1}{2a^2+2a+1} \exp {\left( 1 + \frac{1}{2a} \right)} .
\end{equation}
The function $g$ is decreasing on $(0, \infty)$. Therefore,
$\phi(R) < g(a) < g(1.8) < 1$ .
\end{case}

\begin{case}
Suppose
$a \le 1.8$. From \eqref{b1} and the obvious inequality $b \le a$ we obtain
\begin{equation}
\phi(R) < h(a) =
\frac{2a^2(a+1)}{(3a+2)(2a^2+2a+1)} \exp{\frac{3}{2}} \, .
\end{equation}
The function $h$ is increasing on $(0, \infty)$. Therefore, $\phi(R) < h(a) < h(1.8)
< 1$ which completes the proof.
\end{case}
\end{proof}

\begin{cor} \label{cor:NouInf2}  
Let $f$ be a monic polynomial of degree $n \ge 2$.
Suppose that $z^0$ is an initial point in $\Cset^n $ satisfying
\[
\left \| \frac{W(z^0)}{d(z^0)} \right \|_{\infty} \le \frac{1}{1.4 n + 2.8} \, .
\]
Then Nourein's method \eqref{Nou} is convergent to a root-vector $\xi$ of $f$ with
the fourth order of convergence. Moreover, the estimates \eqref{priN} and
\eqref{posN} hold with $p=\infty$.
\end{cor}

\begin{proof}
Define $\phi_n = \phi(1/(1.4 n + 2.8))$, where $\phi$ is defined by \eqref{phiN} with
$p=\infty$. The sequence $(\phi_n )$ is increasing for $2 \le n \le 19$ and
decreasing for $n \ge 19$. Hence $\phi_n  \le \phi_{19} < 1$ for $n \ge 2$. Now the
conclusion follows from Theorem~\ref{thm:Nou}.
\end{proof}

Corollary~\ref{cor:NouInf2} improves Theorem 3 of \textsc{Nedi\'c} \cite{Ned01}.
Note that Nedi\'c has proved that Nourein's method is convergent under the condition
\(
\| W(z^0) \|_{\infty} \le C(n) \delta(z^0)
\)
where
\[
C(n) = \left \{ \begin{array}{ll}
{{1 / (1.64 \, n + 1.944)}} & \text{for } 3 \le n \le 23;\\
{{1 / (1.42 \, n + 8.7)}} & \text{for } n > 23 .
\end{array} \right.
\]

\begin{cor} \label{cor:NouL1}  
Let $f$ be a monic polynomial of degree $n \ge 2$.
Suppose that $z^0$ is an initial point in $\Cset^n $ satisfying
\[
\left \| \frac{W(z^0)}{d(z^0)} \right \|_{1} \le R = 0.2825 \ldots ,
\]
where $R$  is the unique solution of the equation
\begin{equation} \label{EqL1N}
\frac{x^3}{(1 - 3 x + x^2)^2} \, \exp{\frac{x - x^2}{1 - 3 x + x^2}} = 1
\end{equation}
in the interval $(0,2 / (3 + \sqrt{5})$.
Then $f$ has only simple zeros and Nourein's sequence \eqref{Nou} converges
to a root-vector $\xi$ of $f$ with the fourth order of convergence
and error bounds \eqref{priN} and \eqref{posN} with $p=1$.
\end{cor}

\begin{proof}
By Theorem~\ref{cor:EhrL1} $f$ has only simple zeros. It is easy to show that
$\phi(x) < g(x)$ for $0 < x < 2 / (3 + \sqrt{5})$ , where $g(x)$ denotes the
left-hand side of the equation \eqref{EqL1N} and $\phi(x)$ is defined by \eqref{phiN}
with $p=1$. Therefore, $\phi(R)< 1$ which according to Theorem~\ref{thm:Nou}
completes the proof.
\end{proof}

\noindent
Faculty of Mathematics and Informatics\\
University of Plovdiv\\
Plovdiv 4000, Bulgaria\\
E-mail: proinov@pu.acad.bg

\end{document}